\documentclass[12pt,reqno]{amsart}

\usepackage{graphicx}
\graphicspath{ {./images/} }
\usepackage{bm}

\usepackage{color}

\usepackage{amssymb}

\usepackage{amsfonts}

\usepackage{amsmath}

\usepackage[all]{xy}

\newtheorem{theorem}{Theorem}[section]

\newtheorem{corollary}[theorem]{Corollary}

\newtheorem{lemma}[theorem]{Lemma}

\newtheorem{Definition}[theorem]{Definition}

\newtheorem{Example}[theorem]{Example}

\newtheorem{Remark}[theorem]{Remark}

\newenvironment{remark}{\begin{Remark}\begin{em}}{\end{em}\end{Remark}}

\DeclareMathOperator{\tr}{tr}

\address{Sejong Kim \\ Department of Mathematics, Chungbuk National University, Cheongju 28644, Korea}
\email{skim@chungbuk.ac.kr}

\address{Vatsalkumar N. Mer \\ Institute for Industrial and Applied Mathematics, Chungbuk National University, Cheongju 28644, Korea}
\email{vnm232657@gmail.com}

\begin{document}

\title{Mean iteration with alternative means}

\author{Sejong Kim and Vatsalkumar N. Mer}

\begin{abstract}
Motivated from mean iteration between two Kubo-Ando operator means, we prove the convergence of mean iteration between the alternative means corresponding to two operator monotone functions. Alternative means are binary operations containing the spectral geometric mean and Wasserstein mean. One interesting consequence is that the common limit of mean iteration between the Wasserstein mean and its adjoint is equal to the spectral geometric mean. We also show the convergence of mean iteration between Kubo-Ando operator mean and alternative mean of positive definite matrices.
\end{abstract}



\date{}
\maketitle

\section{Introduction}

Kubo and Ando established in \cite{KA} an axiomatic theory of operator means of positive definite operators. Their remarkable consequence is that an operator mean $\sigma$ is characterized by an operator monotone function, in the sense that
\begin{displaymath}
A \sigma B = A^{1/2} f(A^{-1/2} B A^{-1/2}) A^{1/2},
\end{displaymath}
where $f$ is the operator monotone function on $(0, \infty)$ with $f(1) = 1$.
We often use the notation $\sigma_{f}$ instead of $f$ to indicate the representing function $f$.
For positive definite operators $A, B$ the classical examples of Kubo-Ando's operator mean are the weighted arithmetic, geometric, and harmonic mean given by
\begin{displaymath}
A \nabla_{t} B = (1-t) A + t B, \quad A \sharp_{t} B = A^{1/2} (A^{-1/2} B A^{-1/2})^{t} A^{1/2}, \quad A !_{t} B = ((1-t) A^{-1} + t B^{-1})^{-1}
\end{displaymath}
for $t \in [0,1]$. The arithmetic-geometric-harmonic mean inequalities with respect to the Loewner order are well-known:
\begin{displaymath}
A !_{t} B \leq A \sharp_{t} B \leq A \nabla_{t} B.
\end{displaymath}

Another consequence of Kubo-Ando's operator mean is the mean iteration, called the composition \cite[Section 6]{KA}.
In other words, two sequences $\{ A_{k} \}$ and $\{ B_{k} \}$ defined as $A_{0} = A, \ B_{0} = B$ and
\begin{displaymath}
A_{k+1} = A_{k} \sigma_{f} B_{k}, \quad B_{k+1} = A_{k} \sigma_{g} B_{k}, \quad k \geq 0
\end{displaymath}
with normalized operator monotone functions $f, g$ on $(0, \infty)$ converge to the same limit.
In particular, by letting $f(x) = (1-t) + tx$ and $g(x) = ((1-t) + tx^{-1})^{-1}$
one can see that $A_{k}, B_{k}$ converge to the same limit, named the Gauss mean.
Indeed, $A_{k}$ is a decreasing sequence bounded below by $B_{1} = A !_{t} B$ and $B_{k}$ is an increasing sequence bounded above by $A_{1} = A \nabla_{t} B$ with respect to the Loewner order:
\begin{displaymath}
B_{k} \leq B_{k+1} \leq A_{k+1} \leq A_{k}, \quad k \geq 1.
\end{displaymath}
Especially, for $t = 1/2$ the common limit is the geometric mean $A \sharp B = A \sharp_{1/2} B$.

Recently, the spectral geometric mean $A \natural_{t} B$ and Wasserstein mean $A \diamond_{t} B$ as examples of non-Kubo-Ando's operator mean have been actively studied in \cite{GH, GK24, GK, HK22, Kim}:
\begin{displaymath}
\begin{split}
A \natural_{t} B & = (A^{-1} \sharp B)^{t} A (A^{-1} \sharp B)^{t}, \\
A \diamond_{t} B & = (1-t)^{2} A + t^{2} B + t (1-t) [A (A^{-1} \sharp B) + (A^{-1} \sharp B) A].
\end{split}
\end{displaymath}
New scheme of operator mean containing spectral geometric mean and Wasserstein mean, named an alternative mean, has been introduced in \cite{DFKC}.
For an operator monotone function $f: (0, \infty) \to (0, \infty)$ with $f(1) = 1$, the alternative mean $A \hat{\sigma}_{f} B$ is defined by
\begin{displaymath}
A \hat{\sigma}_{f} B := f(A^{-1} \sharp B) A f(A^{-1} \sharp B).
\end{displaymath}
For given $t \in [0,1]$ the alternative mean corresponding to $f(x) = x^t$ is the spectral geometric mean $A \natural_t B$,
and the alternative mean corresponding to $f(x)= (1-t) + t x$ is the Wasserstein mean $A \diamond_t B$. See \cite{DFK, DFK-2, DFM} for more results.

From the Kubo-Ando's result on mean iteration, it arises a natural question whether or not similar iteration constructed by alternative means corresponding to operator monotone functions $f(x)$ and $g(x)$ converge.
We prove the convergence of this iteration, and particularly observe that the common limit of mean iteration between the Wasserstein mean and its adjoint for $t = 1/2$ coincides with the spectral geometric mean.
Finally, we show the convergence of mean iteration between the Kubo-Ando's operator mean and alternative mean given by the common operator monotone function.

\section{Alternative means}

Let $B(\mathcal{H})$ be the Banach space of all bounded linear operators
on a Hilbert space $\mathcal{H}$ with inner product $\langle \cdot, \cdot \rangle$.
We denote as $S(\mathcal{H})$ the subset of $B(\mathcal{H})$ consisting of all self-adjoint operators.
We say that $A \in S(\mathcal{H})$ is positive semi-definite (definite) if $\langle x, Ax \rangle \geq 0$ ($ > 0$) for any (non-zero, resp.) $x \in \mathcal{H}$.
Let $\mathbb{P}$ be the set of all positive definite invertible operators.
For $X, Y \in S(\mathcal{H})$, we write $X \leq Y$ if and only if $Y-X \geq 0$, and $X < Y$ if and only if $Y-X > 0$.
This relation is known as the Loewner partial order.

In this section we recall the alternative mean with its properties.
For given operator monotone function $f: (0, \infty) \to (0, \infty)$ with $f(1) = 1$,
we define the alternative mean $\hat{\sigma}_{f}$ on $\mathbb{P}$ \cite{DFKC} by
\begin{displaymath}
A \hat{\sigma}_{f} B := f(A^{-1} \sharp B) A f(A^{-1} \sharp B).
\end{displaymath}
By the Riccati identity, one can see that $A \hat{\sigma}_{f} B$ is a unique positive definite solution to the equation
\begin{displaymath}
A^{-1} \sharp X = f(A^{-1} \sharp B).
\end{displaymath}
Alternative mean does not hold the monotonicity with respect to the Loewner order, meanwhile the in-betweenness with respect to the near-order $\preceq$ is satisfied \cite[Theorem 15]{DFKC}:
\begin{center}
$A \preceq A \hat{\sigma}_{f} B \preceq B$ \quad whenever \quad $A \preceq B$,
\end{center}
where $A \preceq B$ if and only if $A^{-1} \sharp B \geq I$.

\begin{theorem}
Let $f, g: (0, \infty) \to (0, \infty)$ be operator monotone functions with $f(1) = g(1) = 1$.
For $A, B \in \mathbb{P}$,
\begin{itemize}
  \item[(1)] $(c A) \hat{\sigma}_{f} (c B) = c (A \hat{\sigma}_{f} B)$ for any $c > 0$,
  \item[(2)] $(U A U^{*}) \hat{\sigma}_{f} (U B U^{*}) = U (A \hat{\sigma}_{f} B) U^{*}$ for any unitary operator $U$,
  \item[(3)] $A \hat{\sigma}_{f} (A \hat{\sigma}_{g} B) = A \hat{\sigma}_{f \circ g} B$.
\end{itemize}
\end{theorem}

In the finite-dimensional setting $\mathcal{H} = \mathbb{C}^{n}$,
let $\mathbb{H}_{n}$ be a real vector space of all Hermitian matrices and
let $\mathbb{P}_{n}$ be the open convex cone of positive definite matrices.
For $t \in [0,1]$, the Wasserstein mean of $A, B \in \mathbb{P}_{n}$, given by
\begin{equation} \label{E:Wass}
\begin{split}
A \diamond_{t} B & = (1-t)^{2} A + t^{2} B + t(1-t) \left[ A (A^{-1} \sharp B) + (A^{-1} \sharp B) A \right] \\
& = \left[ I \nabla_{t} (A^{-1} \sharp B) \right] A \left[ I \nabla_{t} (A^{-1} \sharp B) \right],
\end{split}
\end{equation}
is one of the examples of alternative mean with $f(x) = (1-t) + tx$.
Moreover, $A \diamond_{t} B$ is a least squares mean with respect to the Bures-Wasserstein distance $d_{W}$:
\begin{displaymath}
A \diamond_{t} B = \underset{X \in \mathbb{P}_{n}}{\arg \min} \, (1-t) d_{W}^{2} (A, X) + t d_{W}^{2} (B, X),
\end{displaymath}
where $d_{W}(A, B) = \left[ \tr (A + B) - 2 \tr (A^{1/2} B A^{1/2})^{1/2} \right]^{1/2}$.
Note from \cite[Theorem 1]{BJL} that the Bures-Wasserstein distance has the following extremal characterization:
\begin{equation} \label{E:extremal}
d_{W}(A, B) = \underset{U \in \mathbb{U}_{n}}{\min} \Vert A^{1/2} - B^{1/2} U \Vert_{2},
\end{equation}
where $\mathbb{U}_{n}$ is the unitary group and the minimum is attained by the unitary part for polar decomposition of $B^{1/2} A^{1/2}$.
By vanishing the gradient of objective function $F(X) := (1-t) d_{W}^{2} (A, X) + t d_{W}^{2} (B, X)$ one can see that $A \diamond_{t} B$ is a unique positive definite solution $X$ to the equation
\begin{equation} \label{E:Wass-eq}
(1-t) (A \sharp X^{-1}) + t (B \sharp X^{-1}) = I.
\end{equation}
Hwang and Kim in \cite[Theorem 5.1]{HK19} have defined the Wasserstein mean $A \diamond_{t} B$ of positive definite operators $A, B$ by showing that \eqref{E:Wass-eq} has a unique solution $X \in \mathbb{P}$.


Differently from the Kubo-Ando's geometric mean, $(A^{-1} \diamond_{t} B^{-1})^{-1} \neq A \diamond_{t} B$ in general.
We call $(A^{-1} \diamond_{t} B^{-1})^{-1}$ an adjoint of Wasserstein mean. On the other hand, there are interesting relationships with respect to the near-order.

\begin{lemma} \cite[Lemma 2.8]{HK22} \label{L:self-dual}
Let $A, B \in \mathbb{P}$ and $t \in [0,1]$. Then
\begin{displaymath}
A \diamond_{t} B = \left[ I \nabla_{s} (X \nabla X^{-1}) \right] (A^{-1} \diamond_{t} B^{-1})^{-1} \left[ I \nabla_{s} (X \nabla X^{-1}) \right],
\end{displaymath}
where $s = 2t(1-t) \in [0,1/2]$ and $X = A^{-1} \sharp B$.
Moreover,
\begin{center}
$A \diamond_{t} B = (A^{-1} \diamond_{t} B^{-1})^{-1}$ \ for $t \in (0,1)$ \quad if and only if \quad $A = B$.
\end{center}
\end{lemma}

\begin{lemma} \cite{BJL,DFKC} \label{L:chain}
Let $f$ be an operator monotone function on $(0, \infty)$ with $f(1) = 1$.
For $A, B \in \mathbb{P}$
\begin{equation} \label{E:chain}
A !_{t} B \leq (A^{-1} \diamond_{t} B^{-1})^{-1} \preceq A \hat{\sigma}_{f} B \preceq A \diamond_{t} B \leq A \nabla_{t} B,
\end{equation}
where $f'(1) = t \in [0,1]$ is the weight of the alternative mean $\hat{\sigma}_{f}$.
\end{lemma}

Lemma \ref{L:chain} generalized the near-order inequality $A \natural_{t} B \preceq A \diamond_{t} B$ obtained by Gan and Huang \cite{GH}.
The near-order does not hold the transitivity \cite{DF} so
\begin{displaymath}
(A^{-1} \diamond_{t} B^{-1})^{-1} \preceq A \diamond_{t} B
\end{displaymath}
can not be derived from \eqref{E:chain} but it has been proved in \cite[Theorem 3.9]{GH}.

Let $\lambda_{1}(X), \dots, \lambda_{n}(X)$ be (real) eigenvalues of $X \in \mathbb{H}_{n}$ in decreasing order: $\lambda_{1}(X) \geq \cdots \geq \lambda_{n}(X)$.
The following expression
\begin{displaymath}
\mu_{p}(a, b) = ((1-t) a^{p} + t b^{p})^{1/p}, \quad p \neq 0
\end{displaymath}
is known as the quasi-arithmetic mean of $a, b > 0$.

\begin{lemma} \label{L:Wass-bounds}
Let $A, B \in \mathbb{P}_{n}$ and $t \in [0,1]$. Then
\begin{displaymath}
\mu_{1/2}(\lambda_{n}(A), \lambda_{n}(B)) I \leq A \diamond_{t} B \leq \mu_{1/2}(\lambda_{1}(A), \lambda_{1}(B)) I.
\end{displaymath}
\end{lemma}

\begin{proof}
Note from \eqref{E:Wass-eq} that $A \diamond_{t} B$ is a unique positive definite solution of the equation
\begin{displaymath}
(1-t) (X^{1/2} A X^{1/2})^{1/2} + t (X^{1/2} B X^{1/2})^{1/2} = X.
\end{displaymath}
Since $A \leq \lambda_{1}(A) I$ and $B \leq \lambda_{1}(B) I$, we have by monotonicity of the square root map
\begin{displaymath}
X \leq \left[ (1-t) \sqrt{\lambda_{1}(A)} + t \sqrt{\lambda_{1}(B)} \right] X^{1/2}. \end{displaymath}
Simplifying for $X$, we obtain $X \leq \mu_{1/2}(\lambda_{1}(A), \lambda_{1}(B)) I$.
Similarly, we can prove the other inequality.
\end{proof}

\begin{remark}
For $A_{1}, \dots, A_{m} \in \mathbb{P}_{n}$ and a positive probability vector $\omega = (w_{1}, \dots, w_{m})$ the Wasserstein mean $\Omega(\omega; A_{1}, \dots, A_{m})$ is a least squares mean with respect to the Bures-Wasserstein distance $d_{W}$:
\begin{displaymath}
\Omega(\omega; A_{1}, \dots, A_{m}) = \underset{X \in \mathbb{P}_{n}}{\arg \min} \sum_{j=1}^{m} w_{j} d_{W}^{2} (A_{j}, X).
\end{displaymath}
It is known that $\Omega(\omega; A_{1}, \dots, A_{m})$ coincides with a unique positive definite solution of the gradient-vanishing equation
\begin{displaymath}
\sum_{j=1}^{m} w_{j} A_{j} \sharp X^{-1} = I, \quad \textrm{equivalently} \quad \sum_{j=1}^{m} w_{j} (X^{1/2} A_{j} X^{1/2})^{1/2} = X.
\end{displaymath}
Using the same proof strategy of Lemma \ref{L:Wass-bounds}, we can obtain
\begin{displaymath}
\mu_{1/2}(\lambda_{n}(A_{1}), \dots, \lambda_{n}(A_{m})) I \leq \Omega(\omega; A_{1}, \dots, A_{m}) \leq \mu_{1/2}(\lambda_{1}(A_{1}), \dots, \lambda_{1}(A_{m})) I.
\end{displaymath}
This is an improvement of \cite[Lemma 2.4]{KL}.
\end{remark}

A bound for the operator norm of difference between the arithmetic and Wasserstein means has been shown in \cite[Lemma 3.11]{HK22}:
\begin{displaymath}
\Vert A \nabla_{t} B - A \diamond_{t} B \Vert \leq 2t(1-t) \frac{\beta (\sqrt{\beta} - \sqrt{\alpha})}{\sqrt{\alpha}}
\end{displaymath}
for $A, B \in \mathbb{P}_{n}$ such that $\alpha I \leq A, B \leq \beta I$ for some positive constants $\alpha, \beta$.
We provide another bound for the Frobenius norm. Here, $\sigma(A)$ denotes the spectrum of $A$.
\begin{lemma} \label{L:A-W distance}
Let $A, B \in \mathbb{P}_{n}$ and $t \in [0,1]$. Then
\begin{displaymath}
\Vert A \nabla_{t} B - A \diamond_{t} B \Vert_{2} \leq \frac{t(1-t)}{4 \alpha} \Vert A - B \Vert_{2}^{2},
\end{displaymath}
where $\alpha := \min \{ \lambda_{n}(A), \lambda_{n}(B) \}$.
\end{lemma}

\begin{proof}
By definition \eqref{E:Wass} of the Wasserstein mean,
\begin{displaymath}
A \nabla_{t} B - A \diamond_{t} B = t(1-t) \left[ A + B - A (A^{-1} \sharp B) - (A^{-1} \sharp B) A \right].
\end{displaymath}
Since $\tr (A (A^{-1} \sharp B)) = \tr (A^{1/2} B A^{1/2})^{1/2} = \tr ((A^{-1} \sharp B) A)$,
\begin{displaymath}
\tr \left( A \nabla_{t} B - A \diamond_{t} B \right) = t(1-t) d_{W}^{2}(A, B)
\end{displaymath}
One can easily observe that
\begin{equation} \label{E:BW-bounds}
d_{W}(A, B) \leq \Vert A^{1/2} - B^{1/2} \Vert_{2} \leq \frac{1}{2 \sqrt{\alpha}} \Vert A - B \Vert_{2},
\end{equation}
where the first inequality holds from \eqref{E:extremal} with $U = I$, and the second inequality follows from the Lipschitz continuity of square root map.
Since $A \diamond_{t} B \leq A \nabla_{t} B$ by Lemma \ref{L:chain}, we have $\Vert A \nabla_{t} B - A \diamond_{t} B \Vert_{1} = \tr \left( A \nabla_{t} B - A \diamond_{t} B \right)$.
Since $\Vert X \Vert_{2} \leq \Vert X \Vert_{1}$ for any matrix $X$, we can conclude from \eqref{E:BW-bounds}
\begin{displaymath}
\Vert A \nabla_{t} B - A \diamond_{t} B \Vert_{2} \leq \tr \left( A \nabla_{t} B - A \diamond_{t} B \right) = t(1-t) d_{W}^{2}(A, B) \leq \frac{t(1-t)}{4 \alpha} \Vert A - B \Vert_{2}^{2}.
\end{displaymath}
\end{proof}

\section{Alternative mean iteration}

Let $f, g: (0, \infty) \to (0, \infty)$ be operator monotone functions with $f(1) = g(1) = 1$.
If operator monotone functions $f$ is said to be trivial if $f(x) = 1$ or $f(x) = x$.
Let $A, B \in \mathbb{P}$.
In this section, we consider two sequences $\{ A_{k} \}_{k \geq 0}$ and $\{ B_{k} \}_{k \geq 0}$ constructed by alternative means,
such as $A_{0} = A, \ B_{0} = B$ and
\begin{equation} \label{E:Alter-sequences}
A_{k+1} = A_{k} \hat{\sigma}_{f} B_{k}, \quad B_{k+1} = A_{k} \hat{\sigma}_{g} B_{k},
\end{equation}
where one of $f$ and $g$ is non-trivial.

For $x > 0$ we define $\displaystyle r(x) := f(x)^{-1} g(x) $.
For $X \in \mathbb{P}$, $r(X) = f(X)^{-1} g(X) = g(X) f(X)^{-1}$, since $f(X)$ and $g(X)$ commute.
We denote as $f^{k}$ the $k$-times decomposition of $f$ with $f^{0}(x) = x$.

\begin{lemma} \cite[Lemma 5.3]{KA} \label{L:KA-Lemma}
Let $f, g: (0, \infty) \to (0, \infty)$ be operator monotone functions with $f(1) = g(1) = 1$. If one of $f$ and $g$ is non-trivial, then for any $x >0$
\begin{displaymath}
\lim_{k \to \infty} r^{k}(x) = 1
\end{displaymath}
in operator norm.
\end{lemma}

\begin{theorem} \label{thm:main}
Let $A, B \in \mathbb{P}$.
Let $\{ A_k \}$ and $\{ B_k \}$ be the sequences given by \eqref{E:Alter-sequences}.
Then
\begin{equation} \label{eq:common-limit}
\lim_{k \to \infty} A_k = \lim_{k \to \infty} B_k.
\end{equation}
Indeed, the common limit is also the alternative mean $L(A, B) := f_{\infty}(T) A f_{\infty}(T)$ with $T := A^{-1} \sharp B$
associated to the operator monotone function $f_{\infty}$ defined by
\begin{displaymath}
f_{\infty}(x) = \lim_{k \to \infty} f_{k}(x) f(r^{k}(x)),
\end{displaymath}
where $f_{k+1}(x) = f_{k}(x) f(r^{k}(x))$ for $k \geq 0$ and $f_{0}(x) = f(x)$.
\end{theorem}

\begin{proof}
Define $T_k := A_k^{-1} \sharp B_k$ for $k \geq 0$.
We have $T_0 = A^{-1} \sharp B = T$. Note that
\begin{displaymath}
\begin{split}
T_1 A_1 T_1 & = T_1 (A \hat{\sigma}_{f} B) T_1 = T_{1} f(T_0) A f(T_0) T_{1}, \\
B_{1} & = g(T_{0}) A g(T_{0}) = f(T_0)^{-1} g(T_0) f(T_0) A f(T_0) g(T_0) f(T_0)^{-1}.
\end{split}
\end{displaymath}
The second equality for $B_{1}$ follows from the commutativity of $f(T_{0})$ and $g(T_{0})$.
Since $T_{1} A_{1} T_{1} = B_{1}$ by the Riccati identity corresponding to $T_{1} = A_{1}^{-1} \sharp B_{1}$,
we obtain
$T_1 = f(T_0)^{-1} g(T_0) = r(T_0)$.
Applying the same argument recursively, we obtain
\begin{equation} \label{eq:T-iteration}
T_{k+1} = r(T_k) = r^{k+1}(T_0).
\end{equation}

We now claim that
\begin{equation} \label{eq:pair-representation}
A_k = f_k (T_0) A f_k(T_0), \qquad B_k = g_k (T_0) A g_k (T_0),
\end{equation}
for $k \geq 1$, where
\begin{displaymath}
f_k(x) := f_{k-1}(x) f(r^{k-1}(x)) \quad \text{and} \quad g_{k}(x) := f_{k-1}(x) g(r^{k-1}(x)), \quad k \geq 1
\end{displaymath}
with $f_{0}(x) = f(x)$ and $g_{0}(x) = g(x)$ for all $x > 0$.
We prove the formula for $A_k$ by induction.
For $k = 1$, we have $A_1 = f_{1}(T_{0}) A f_{1}(T_{0})$ and $B_1 = g_{1}(T_{0}) A g_{1}(T_{0})$ by definition.
If $A_k = f_{k}(T_0) A f_{k}(T_0)$, then
\begin{displaymath}
A_{k+1} = f(T_k) A_{k} f(T_k) = f(T_k) f_{k}(T_0) A f_{k}(T_0) f(T_k)
= f_{k}(T_0) f(r^{k}(T_0)) A f_{k}(T_0) f(r^{k}(T_0)),
\end{displaymath}
where the last equality follows from \eqref{eq:T-iteration} and the commutativity of $f_{k}$ and $f \circ r^{k}$.
Thus, $A_{k+1} = f_{k+1}(T_0) A f_{k+1}(T_0)$.
Since $T_{k} f(T_{k-1}) A_{k-1} f(T_{k-1}) T_{k} = T_k A_k T_k = B_{k} = g(T_{k-1}) A_{k-1} g(T_{k-1})$ by the Riccati identity for $T_{k} = A_{k}^{-1} \sharp B_{k}$, we have
\begin{displaymath}
T_k = (f(r^{k-1}(T_0)))^{-1} g(r^{k-1}(T_0)) = f_k (T_0)^{-1} g_k (T_0),
\end{displaymath}
which implies $g_k (T_0) = f_{k-1}(T_0) g(r^{k-1}(T_0))$ by definition of $f_{k}$.
Therefore, we obtain from the formula for $A_{k}$
\begin{displaymath}
B_k = T_{k} A_{k} T_{k} = f_k (T_0)^{-1} g_k (T_0) f_k (T_0) A f_k(T_0) f_k (T_0)^{-1} g_k (T_0) = g_k (T_0) A g_k (T_0),
\end{displaymath}
which proves the claim.

By the spectral theorem, we write
\begin{displaymath}
T_{0} = \int_{\sigma(T_{0})} \lambda \, dE(\lambda),
\end{displaymath}
where $E$ is the projection-valued spectral measure of $T_{0}$.
Then \eqref{eq:T-iteration} implies $\displaystyle T_{k} = \int_{\sigma(T_{0})} r^{k}(\lambda) \, dE(\lambda)$ by functional calculus, and
\begin{displaymath}
\Vert T_{k} - I \Vert = \underset{x \in \sigma(T_{0})}{\sup} |r^{k}(x) - 1| \to 0 \quad \textrm{as} \quad k \to \infty
\end{displaymath}
by Lemma \ref{L:KA-Lemma}. Thus,
\begin{equation} \label{eq:T-to-I}
T_k \longrightarrow I \quad \text{as} \quad k \to \infty
\end{equation}
in operator norm.
As shown in \cite[Theorem 6.2]{KA} with $a_{k} = f_{k}(x)$ and $b_{k} = g_{k}(x)$, these sequences converge to the same limit.
Therefore, the representation \eqref{eq:pair-representation} gives
\[
A_k = f_k(T_0) A f_k(T_0) \longrightarrow f_\infty (T_0) A f_\infty (T_0) = L(A,B),
\]
and similarly,
\[
B_k = g_k(T_0) A g_k(T_0) \longrightarrow L(A,B).
\]
\end{proof}

\begin{corollary}
Let $A, B \in \mathbb{P}$ and $t \in (0,1)$. Two sequences $\{ A_k \}$ and $\{ B_k \}$ given by $A_{0} = A, \ B_{0} = B$ and
\begin{displaymath}
A_{k+1} = A_{k} \diamond_{t} B_{k}, \quad B_{k+1} = (A_{k}^{-1} \diamond_{t} B_{k}^{-1})^{-1}.
\end{displaymath}
Then $A_{k}, B_{k}$ converge to the same limit as $k \to \infty$.
In particular, for $t = 1/2$
\begin{equation} \label{eq:spectral-limit}
\lim_{k \to \infty} A_k = \lim_{k \to \infty} B_k = A \natural B.
\end{equation}
\end{corollary}

\begin{proof}
Let $f(x) = (1-t) + t x$ for $x > 0$. Then $f$ is a non-trivial operator monotone function on $(0, \infty)$ with $f(1) = 1$, and
\begin{displaymath}
A_{k+1} = A_{k} \hat{\sigma}_{f} B_{k}, \quad B_{k+1} = A_{k} \hat{\sigma}_{g} B_{k},
\end{displaymath}
where $g(x) = f(x^{-1})^{-1}$. By Theorem \ref{thm:main}, the first assertion holds.

When $t = 1/2$, we have
\begin{displaymath}
r(x) = \frac{1}{f(x) f(x^{-1})} = \frac{4x}{(1+x)^2} = \frac{x}{\bigl( (1+x)/2 \bigr)^2}.
\end{displaymath}
If $x_{j+1} = r(x_j)$ for $j \geq 0$, then $\displaystyle \left( \frac{1 + x_j}{2} \right)^2 = \frac{x_j}{x_{j+1}}$.
Consequently, the finite product with $x_{0} := x$ for $x > 0$ reduces to
\[
f_{k}(x)^2
= \prod_{j=0}^{k-1} \left( \frac{1 + x_j}{2} \right)^2
= \frac{x_0}{x_k}.
\]
Since $x_k \to 1$ by Lemma \ref{L:KA-Lemma}, it gives $f_{\infty}(x) = \sqrt{x}$.
Therefore, formula \eqref{eq:spectral-limit} follows from Theorem~\ref{thm:main}.
\end{proof}

\section{Mean iterations with Kubo-Ando and alternative means}

Let us first consider two sequences $\{ A_{k} \}_{k \geq 0}$ and $\{ B_{k} \}_{k \geq 0}$ constructed by the weighted arithmetic mean and Wasserstein mean: $A_{0} = A, B_{0} = B$ and
\begin{equation} \label{E:A-W-sequences}
A_{k+1} = A_{k} \nabla_{t} B_{k}, \quad B_{k+1} = A_{k} \diamond_{t} B_{k},
\end{equation}
where $t \in (0,1)$. From Lemma \ref{L:chain} we obtain
\begin{equation} \label{E:A-W}
B_{k} \leq A_{k}, \quad k \geq 1.
\end{equation}
Since $A_{k+1} \leq A_{k}$ for all $k \geq 1$, $A_{k}$ is a decreasing sequence with respect to the Loewner order and is bounded below by $0$
so it converges, say to $L$ as $k \to \infty$.
Due to \eqref{E:A-W},
$B_{k+1} \succeq B_{k}$ for all $k \geq 1$ from \cite[Theorem 3.2 (6)]{DFK}.
On the other hand, the near-order does not hold the transitivity so we could not say that $B_{k}$ is an increasing sequence.

\begin{remark} \label{R:A-W convergence}
We can prove that $A_{k}$ and $B_{k}$ converge to the same limit as $k \to \infty$ under certain condition on eigenvalues of $A$ and $B$.
Let $\alpha := \min \{ \lambda_{n}(A), \lambda_{n}(B) \}$ and $\beta := \max \{ \lambda_{1}(A), \lambda_{1}(B) \}$.
We first claim that
\begin{center}
$\min \{ \lambda_{n}(A_{k}), \lambda_{n}(B_{k}) \} \geq \alpha$ \ and \ $\max \{ \lambda_{1}(A_{k}), \lambda_{1}(B_{k}) \} \leq \beta$ \ for all $k \geq 0$.
\end{center}
From monotonicity of the weighted arithmetic mean and Lemma \ref{L:Wass-bounds},
\begin{displaymath}
\begin{split}
\alpha \leq [\lambda_{n}(A) \nabla_{t} \lambda_{n}(B)] I_{n} \leq & A \nabla_{t} B \leq [\lambda_{1}(A) \nabla_{t} \lambda_{1}(B)] I_{n} \leq \beta, \\
\alpha \leq \mu_{1/2}(\lambda_{n}(A), \lambda_{n}(B)) I_{n} \leq & A \diamond_{t} B \leq \mu_{1/2}(\lambda_{1}(A), \lambda_{1}(B)) I_{n} \leq \beta.
\end{split}
\end{displaymath}
By using the previous inequalities inductively, we can prove the claim.

From Lemma \ref{L:A-W distance}
\begin{equation} \label{E:norm}
\Vert A \nabla_{t} B - A \diamond_{t} B \Vert_{2} \leq \frac{t(1-t)}{4 \alpha} \Vert A - B \Vert_{2}^{2} \leq \frac{\sqrt{n} (\beta - \alpha)}{16 \alpha} \Vert A - B \Vert_{2},
\end{equation}
since $\displaystyle \Vert A - B \Vert_{2} \leq \left[ \sum_{i=1}^{n} \lambda_{i}(A - B)^{2} \right]^{1/2} \leq \sqrt{n} (\beta - \alpha)$ and $t(1-t) \leq 1/4$.
Using the preceding argument with \eqref{E:norm} inductively,
\begin{displaymath}
\Vert A_{k+1} - B_{k+1} \Vert_{2} \leq \tau^{k+1} \Vert A - B \Vert_{2}, \quad \tau := \frac{\sqrt{n} (\beta - \alpha)}{16 \alpha}.
\end{displaymath}
If $\beta / \alpha < 1 + 16/\sqrt{n}$, then $0 < \tau < 1$ so $A_{k+1} - B_{k+1}$ converges to $0$ as $k \to \infty$.
Note that $B_{k+1} = A_{k+1} - (A_{k+1} - B_{k+1})$ for all $k$.
We have seen that $A_{k+1}$ converges to $L$ as $k \to \infty$. Hence,
\begin{displaymath}
\lim_{k \to \infty} B_{k+1} = \lim_{k \to \infty} A_{k+1} - \lim_{k \to \infty} (A_{k+1} - B_{k+1}) = L.
\end{displaymath}
\end{remark}

Note that $A_{k}, B_{k}$ in \eqref{E:A-W-sequences} are special examples of the Kubo-Ando's operator mean and alternative mean with the operator monotone function $f(x) = (1-t) + t x$ for $x > 0$.
Next, we show the convergence of $A_{k}$ and $B_{k}$ constructed by any operator monotone function $f: (0, \infty) \to (0, \infty)$ with $f(1) = 1$.
All norms below are operator norms, meanwhile $\| \cdot \|_2$ denotes the Frobenius norm.

\begin{lemma} \label{L:trace-ineq}
For $A, B \in \mathbb{P}_{n}$,
\begin{displaymath}
\tr (A \hat{\sigma}_f B) \leq (1-t) \tr A + t \tr B - t(1-t) \tr (A (X - I)^{2}),
\end{displaymath}
where $X = A^{-1} \sharp B$.
\end{lemma}

\begin{proof}
Put $t = f'(1)$. The standard integral representation of normalized operator monotone function $f$ is
\[
 f(x) = \int_{[0,1]} \frac{x}{s + (1-s) x} \, d\mu(s),
 \qquad t = \int_{[0,1]} s \, d\mu(s),
\]
for a probability measure $\mu$.
In particular, $0 \leq t \leq 1$.
The weighted arithmetic-harmonic mean inequality and Jensen's inequality in the variable $s$ give
\begin{equation} \label{eq:scalar-bounds}
 \frac{1}{(1-t)+t/x} \leq f(x) \leq (1-t)+tx \qquad ( x > 0 ).
\end{equation}
Since all quantities in \eqref{eq:scalar-bounds} are positive,
\begin{displaymath}
f(x)^{2} \leq ((1-t) + t x)^{2} = (1-t) + t x^{2} - t(1-t) (x-1)^{2}.
\end{displaymath}
So letting $X = A^{-1} \sharp B$ implies
\begin{displaymath}
\tr (A \hat{\sigma}_f B) = \tr (A f(X)^{2})
\leq (1-t) \tr A + t \tr (A X^{2}) - t(1-t) \tr (A (X - I)^{2}).
\end{displaymath}
Since $\tr (A X^{2}) = \tr (X A X) = \tr B$ by the Riccati identity, it completes the proof.
\end{proof}

\begin{lemma}
Let $P \in \mathbb{P}_{n}$ such that $m I \leq P \leq M I$ for some $m, M > 0$.
For small enough $E \in \mathbb{H}_{n}$ that $P + E > 0$,
\begin{align}
 P \sigma_f (P+E) &= P + t E + O(\| E \|^2), \label{eq:ka-expansion} \\
 P \hat{\sigma}_f (P+E) &= P + t E + O(\| E \|^2). \label{eq:alt-expansion}
\end{align}
\end{lemma}

\begin{proof}
For \eqref{eq:ka-expansion}, set $H := P^{-1/2} E P^{-1/2}$. Since $f(1) = 1$ and $f'(1) = t$, we have from Taylor expansion
\[
 f(I + H) = I + t H + O(\|H\|^2).
\]
The remainder is uniform because $\| P^{-1/2} \| \leq m^{-1/2}$ and $\| P \| \leq M$.

For \eqref{eq:alt-expansion}, write $P^{-1} \sharp (P+E) = I + Z$.
The Riccati identity gives
\begin{equation}\label{eq:sylvester}
 ZP + PZ + ZPZ = E.
\end{equation}
The map $(P,Q) \mapsto P^{-1} \sharp Q$ is smooth on $\mathbb P_n \times \mathbb P_n$. Compactness therefore gives the uniform estimate $Z=O(\|E\|)$ near the diagonal. Equivalently, this estimate follows from the implicit function theorem applied to \eqref{eq:sylvester}: its derivative in $Z$ at $Z=0$ is the invertible Sylvester map $Z \mapsto ZP + PZ$.

Now since $f(I + Z) = I + t Z + O(\|Z\|^2)$,
\begin{align*}
 P \hat{\sigma}_f(P+E) &= f(I + Z) P f(I + Z) \\
 &= P + t (ZP + PZ) + O(\|Z\|^2) \\
 &= P + t E + O(\|E\|^2).
\end{align*}
The last equality follows since $ZP + PZ = E - ZPZ$ from \eqref{eq:sylvester} and $\Vert Z P Z \Vert \leq M \Vert Z \Vert = O(\|E\|^2)$.
All remainder constants are uniform on the compact set $[m, M] := \{ X \in \mathbb{H}_{n}: m I \leq X \leq M I \}$ of $\mathbb{P}_{n}$. The scalar Taylor expansions used here are valid because operator monotone functions are real analytic on $(0,\infty)$; functional calculus transfers their local remainder estimates to Hermitian matrices.
\end{proof}

\begin{theorem} \label{T:main}
For $A, B \in \mathbb P_n$, we define $A_0 = A, B_0 = B$, and
\begin{equation} \label{E:KA-Alter}
 A_{k+1} = A_k \sigma_f B_k, \qquad B_{k+1} = A_k \hat{\sigma}_f B_k
\end{equation}
for $k \geq 0$. Then
\[
 \lim_{k \to \infty} A_k = \lim_{k \to \infty} B_k
\]
in norm.
\end{theorem}

\begin{proof}
Let $t = f'(1)$.
If $t=0$, the representation implies $f \equiv 1$, so $A_1 = B_1 = A$.
If $t=1$, it implies $f(x) = x$, so $A_1 = B_1 = B$.
Both sequences are constant thereafter. Hence, we assume $0<t<1$, without loss of generality.

\noindent \textbf{Step 1.}
For $P, Q \in \mathbb P_n$, $P !_{t} Q \leq P \sigma_f Q \leq P \nabla_{t} Q$ by \eqref{eq:scalar-bounds}.
So
\begin{align}
 \tr(P \sigma_f Q) & \leq (1-t) \tr P + t \tr Q, \label{eq:KA-trace} \\
 \tr \bigl( (P \sigma_f Q)^{-1} \bigr) & \leq (1-t) \tr P^{-1} + t \tr Q^{-1}. \label{eq:KA-inverse}
\end{align}
By \cite[Theorem 4]{DFKC}
\begin{equation} \label{eq:alter-chain}
P !_{t} Q \leq (P^{-1} \diamond_{t} Q^{-1})^{-1} \preceq P \hat{\sigma}_f Q \preceq P \diamond_{t} Q \leq P \nabla_{t} Q.
\end{equation}
Then
\begin{align}
 \tr(P \hat{\sigma}_f Q) & \leq (1-t) \tr P + t \tr Q, \label{eq:Alter-trace} \\
 \tr \bigl( (P \hat{\sigma}_f Q)^{-1} \bigr) & \leq (1-t) \tr P^{-1} + t \tr Q^{-1}. \label{eq:Alter-inverse}
\end{align}
Define
\[
 S_k := (1-t) \tr A_k + t \tr B_k,
 \qquad
 T_k := (1-t) \tr A_k^{-1} + t \tr B_k^{-1}.
\]
Equations \eqref{eq:KA-trace}--\eqref{eq:KA-inverse} and \eqref{eq:Alter-trace}--\eqref{eq:Alter-inverse} imply $S_{k+1} \leq S_k$ and $T_{k+1} \leq T_k$.
Let $a = \min \{ t, 1-t \} > 0$, and set
\[
 m = \frac{a}{T_0}, \qquad M = \frac{S_0}{a},
\]
Since $a \tr A_{k} \leq (1-t) \tr A_{k} \leq S_{k} \leq S_{0}$ for all $k \geq 0$,
we have $\lambda_{\max} (A_{k}) \leq \tr A_{k} \leq M$ so
\begin{displaymath}
A_{k} \leq M I.
\end{displaymath}
On the other hand, since $\lambda_{\min}(A_{k})^{-1} = \lambda_{\max}(A_{k}^{-1}) \leq \tr (A_{k}^{-1}) \leq T_{0}/a$ we get $\lambda_{\min}(A_{k}) \geq m$ for all $k \geq 0$ so
\begin{displaymath}
A_{k} \geq m I.
\end{displaymath}
Similarly, we can obtain $m I \leq B_{k} \leq M I$ for all $k \geq 0$.
This means that all iterates lie in a compact subset $[m, M]$ of $\mathbb P_n$.

\noindent \textbf{Step 2.}
Write $X_k := A_k^{-1} \sharp B_k$ for each $k$.
By \eqref{eq:KA-trace} and Lemma \ref{L:trace-ineq}
\begin{displaymath}
\begin{split}
S_{k+1} & \leq \tr A_{k+1} + t \tr B_{k+1} \\
& \leq (1-t) [(1-t) \tr A_{k} + t \tr B_{k}] + t \left[ (1-t) \tr A_{k} + t \tr B_{k} - t(1-t) \tr \bigl( A_{k} (X_{k} - I)^{2} \bigr) \right] \\
& = (1-t) \tr A_{k} + t \tr B_{k} - t^{2} (1-t) \tr \bigl( A_{k} (X_{k} - I)^{2} \bigr),
\end{split}
\end{displaymath}
which implies
\[
 S_k - S_{k+1}
 \geq t^2(1-t) \tr \bigl( A_k (X_k - I)^2 \bigr).
\]
Then the partial sum gives
\begin{displaymath}
t^2(1-t) \sum_{k=0}^{N} \tr \bigl( A_k (X_k - I)^2 \bigr) \leq \sum_{k=0}^{N} (S_{k} - S_{k+1}) = S_{0} - S_{N+1} \leq S_{0}
\end{displaymath}
since $S_{k} \geq 0$ for all $k$, and
\[
 \sum_{k=0}^{N} \tr \bigl( A_k (X_k - I)^2 \bigr) \leq \frac{S_{0}}{t^{2}(1-t)}.
\]
So the partial sum is increasing and bounded above by $\frac{S_{0}}{t^{2}(1-t)}$, and thus, it converges, that is,
\[
 \sum_{k=0}^{\infty} \tr \bigl( A_k (X_k - I)^2 \bigr) < \infty.
\]
This yields that the sequence $\tr \bigl( A_k (X_k - I)^2 \bigr)$ converges to zero as $k \to \infty$.
Moreover, since $(X_k - I) A_k (X_k - I) \geq m (X_k - I)^{2} \geq 0$,
\[
 \tr \bigl( A_k (X_k - I)^2 \bigr)
 \geq m \| X_k - I \|_{2}^2.
\]
Hence, $X_k \to I$ as $k \to \infty$.
Using $B_k = X_k A_k X_k$ and \textbf{Step 1}, we obtain
\[
 B_k - A_k = (X_k - I) A_k X_k + A_k (X_k - I) \to 0.
\]
Thus,
\begin{equation}\label{eq:gap-zero}
 d_k := \| A_k - B_k \| \longrightarrow 0 \quad \text{as} \quad k \to \infty..
\end{equation}

\noindent \textbf{Step 3.}
Applying \eqref{eq:ka-expansion} and \eqref{eq:alt-expansion} with $P = A_k$ and $E = B_k - A_k (=: E_{k})$ and using \eqref{eq:gap-zero},
we obtain positive constants $C_{1}, C_{2}$ such that
\begin{displaymath}
A_{k+1} = A_{k} + t E_{k} + R_{k}, \quad B_{k+1} = A_{k} + t E_{k} + \hat{R}_{k},
\end{displaymath}
and $\Vert R_{k} \Vert \leq C_{1} d_{k}^{2}$ and $\Vert \hat{R}_{k} \Vert \leq C_{2} d_{k}^{2}$ for sufficiently large $k$.
So $A_{k+1} - B_{k+1} = R_{k} - \hat{R}_{k}$, which implies
\begin{equation} \label{eq:quadratic}
 d_{k+1} \leq C d_k^2,
\end{equation}
where $C = C_{1} + C_{2} (>0)$.
Moreover, there exists a positive constant $C'$ such that
\begin{equation}\label{eq:difference}
 \| A_{k+1} - A_k \| \leq t d_{k} + C_{1} d_{k}^{2} \leq C' d_{k}
\end{equation}
for sufficiently large $k$, because $d_{k} \to 0$ as $k \to \infty$.

Since $d_{k} \to 0$ as $k \to \infty$, we can choose $k_0$ sufficiently large that $C d_k \leq 1/2$ for all $k \geq k_0$. Then $d_{k+1} \leq \frac{1}{2} d_k$ for all $k \geq k_{0}$ so
\begin{displaymath}
 \sum_{k=k_0}^{\infty} d_k = k_{0} \sum_{k=0}^{\infty} \frac{1}{2^{k}} < \infty.
\end{displaymath}
The second estimate in \eqref{eq:quadratic} shows that $A_k$ is Cauchy so it converges, say to $L$.
By \eqref{eq:gap-zero}, $B_k \to L$ as well.
\end{proof}

We investigate the properties of the common limit of $A_{k}$ and $B_{k}$ appeared in Theorem \ref{T:main}.
We denote it as $\widehat{L}_{f}(A, B)$.
\begin{theorem} \label{T:properties}
Let $f: (0, \infty) \to (0, \infty)$ be an operator monotone function with $f(1) = 1$.
Let $A, B \in \mathbb P_n$.
\begin{itemize}
  \item[(1)] $\widehat{L}_{f}(A, A) = A$.
  \item[(2)] $\widehat{L}_{f}(c A, c B) = c \widehat{L}_{f}(A, B)$ for any $c > 0$.
  \item[(3)] $U \widehat{L}_{f}(A, B) U^{*} = \widehat{L}_{f}(U A U^{*}, U B U^{*})$ for any unitary matrix $U$.
  \item[(4)] $\widehat{L}_{f}(A, B) = \widehat{L}_{f}(B, A)$ if $f(x) = x f(x^{-1})$ for all $x > 0$.
  \item[(5)] $\widehat{L}_{f}(A, B) = \widehat{L}_{f}(A^{-1}, B^{-1})^{-1}$ if $f(x) = f(x^{-1})^{-1}$ for all $x > 0$.
\end{itemize}
\end{theorem}

\begin{proof}
Idempotency (1), homogeneity (2), and unitary congruence invariance (3) are obvious, because both the Kubo-Ando operator mean and alternative mean satisfy commonly.

By \cite[Corollary 4.2]{KA} and \cite[Theorem 7]{DFKC}
\begin{displaymath}
A \sigma_{f} B = B \sigma_{f} A \quad \textrm{and} \quad A \hat{\sigma}_{f} B = B \hat{\sigma}_{f} A
\end{displaymath}
if and only if $f(x) = x f(x^{-1})$ for all $x > 0$.
Moreover, by \cite[Corollary 4.3]{KA} and \cite[Theorem 6]{DFKC}
\begin{displaymath}
A \sigma_{f} B = (A^{-1} \sigma_{f} B^{-1})^{-1} \quad \textrm{and} \quad A \hat{\sigma}_{f} B = (A^{-1} \hat{\sigma}_{f} B^{-1})^{-1}
\end{displaymath}
if and only if $f(x) = f(x^{-1})^{-1}$ for all $x > 0$.
These yield the symmetry (4) and self-duality (5) for $\widehat{L}_{f}$.
\end{proof}

\begin{remark}
The convergence of $A_{k}$ and $B_{k}$ given by \eqref{E:A-W-sequences} has been shown in Remark \ref{R:A-W convergence} under the condition of spectrum of $A$ and $B$, while Theorem \ref{T:main} provides it without any assumption.
One may ask whether the common limit of $A_{k}$ and $B_{k}$ satisfies the monotonicity for certain order relation.
Indeed, let
\[
A = \begin{pmatrix} 1 & -1 \\ -1 & 2 \end{pmatrix}, \qquad
B = \begin{pmatrix} 277 & 41 \\ 41 & 10 \end{pmatrix}, \qquad
C = \begin{pmatrix} 278 & 30 \\ 30 & 131 \end{pmatrix}.
\]
These matrices are positive definite, and $B \leq C$.
Numerical computation to find the common limit $\widehat{L}_{f}$ in \eqref{E:A-W-sequences} with $f(x) = (1+x)/2$ gives
\begin{displaymath}
\begin{split}
\widehat{L}_{f}(A,B) & \approx
\begin{pmatrix}
107.144826213182 & 13.468434730728 \\
13.468434730728 & 4.591510084501
\end{pmatrix}, \\
\widehat{L}_{f}(A,C) & \approx
\begin{pmatrix}
107.200490580438 & 9.254066966120 \\
9.254066966120 & 53.047954791877
\end{pmatrix}.
\end{split}
\end{displaymath}
Their eigenvalues in decreasing order are approximately
\begin{displaymath}
\begin{split}
\lambda(\widehat{L}_{f}(A,B))
& = (108.884150610201, \; 2.852185687482) \\
\lambda(\widehat{L}_{f}(A,C))
& = (108.738240902694, \; 51.510204469621).
\end{split}
\end{displaymath}
In particular, $\lambda_1(\widehat{L}_{f}(A,B)) > \lambda_1(\widehat{L}_{f}(A,C))$,
which yields
\[
B \leq C \qquad \text{but} \qquad
\Vert \widehat{L}_{f}(A,B) \Vert \not\leq \Vert \widehat{L}_{f}(A,C) \Vert.
\]
This implies that any inequality between $\widehat{L}_{f}(A,B)$ and $\widehat{L}_{f}(A,C)$
with respect to the Loewner order and near order fails.
The computations were performed with $70$-digit working precision;
for both initial pairs, the Frobenius norm of the difference
between the two iterates was less than $10^{-60}$ after seven
iterations.
\end{remark}

\begin{remark}
Let us consider sequences $A_{k}$ and $B_{k}$ given by $A_{0} = A, B_{0} = B$ and
\begin{equation} \label{E:geomean-iteration}
A_{k+1} = A_{k} \sharp_{t} B_{k}, \quad B_{k+1} = A_{k} \natural_{t} B_{k} \quad (k \geq 0),
\end{equation}
where $t \in (0,1)$. Since their common representing function is $f(x) = x^{t}$, Theorem \ref{T:main} yields that they converge to the same limit.
From the log-majorization relationship between the Kubo-Ando geometric mean and spectral geometric mean:
\begin{displaymath}
A \sharp_{t} B \prec_{\log} \exp [(1-t) \log A + t \log B] \prec_{\log} A \natural_{t} B,
\end{displaymath}
one can naturally ask whether or not the common limit is equal to the log-Euclidean mean $\mathrm{LE}(t; A, B) := \exp [(1-t) \log A + t \log B]$.
Indeed, let
\[
A = \begin{pmatrix} 1 & -1 \\ -1 & 2 \end{pmatrix}, \qquad
B = \begin{pmatrix} 2 & 3 \\ 3 & 5 \end{pmatrix}.
\]
Then $A, B$ are positive definite matrices with determinant $1$.
Numerical computation for iteration \eqref{E:geomean-iteration} with $t = 1/2$ gives
\begin{displaymath}
\widehat{L}_{f}(A,B) \approx
\begin{pmatrix}
0.62708672 & 0.48375519 \\
0.48375519 & 1.96786031
\end{pmatrix}
\ \text{but} \
\mathrm{LE}(t; A, B) \approx
\begin{pmatrix}
0.60406852 & 0.47491716 \\
0.47491716 & 2.02882000
\end{pmatrix}.
\end{displaymath}
These means are different, even though their determinant are equal to $1$.
This inspires from Theorem \ref{T:properties} that mean iteration between Kubo-Ando operator mean and alternative mean gives another type of matrix mean.
\end{remark}

We close our paper with an open question.
Alternative mean iteration has been shown in the setting $\mathbb{P}$, 
but iteration between Kubo-Ando mean and alternative mean is proved on the finite-dimensional setting $\mathbb{P}_{n}$.
So a natural question arises whether Theorem \ref{T:main} holds on $\mathbb{P}$.

\section*{Acknowledgments}
The work of Vatsalkumar N. Mer  was supported by Basic Science Research Program through the National Research Foundation of Korea (NRF) funded by the Ministry of Education, Korea (No. RS-2024-00462498).

\end{document}